\documentclass[a4paper,11pt]{article}
\usepackage{amssymb}
\usepackage{latexsym}

\newtheorem{defi}{Definition}

\newtheorem{prop}{Proposition}
\newtheorem{lemma}{Lemma}
\newtheorem{cor}{Corollary}
\newtheorem{lem}{Lemma}[section]

\title{{\sc Fedosov Star-Products} and {\sc 1-Differentiable Deformations}}

\author{Philippe Bonneau
           \thanks{\scriptsize D\'epartement de Math\'ematiques,
Universit\'e de Bourgogne, BP 400, F-21011 Dijon Cedex, France.
            E-mail: {\tt bonneau@u-bourgogne.fr}  }  }

\date{}

\begin{document}

\maketitle

\begin{abstract}
We show that every star product on a symplectic manifold defines uniquely
a 1-differentiable deformation of the Poisson bracket. Explicit formulas
are given. As a corollary we can identify the characteristic class of any
star product as a part of its explicit (Fedosov) expression.
\end{abstract}



\section{Introduction}

The 1-differentiable deformations of the Poisson bracket Lie algebra of
differentiable functions on a Poisson manifold $M$ are usual formal
deformations (in the sense of \cite{G}) built using exclusively
(1,1)-bidifferential operators as cochains. They define a formal Poisson
structure on $M$ starting with the initial Poisson structure
(contravariant 2-tensor) of $M$. In 1974, in one of the first papers of what
 became the {\it deformation quantization theory} (for a review see \cite{W}),
 M. Flato, A. Lichnerowicz and D. Sternheimer \cite{FLS} studied these
deformations for a symplectic manifold $M$.
They showed in particular that the infinitesimal 1-differentiable deformations
(with cochains vanishing on constants) are exactly classified by the second
de Rham cohomology space of $M$.\\
\noindent
As is now well-known, there is a similar classification for star products,
given by a sequence of de Rham 2-cocycles. Comparing both results one
suspects that the ``difference" (the ``difference of what is added on Moyal")
between two star products on $M$ is made of a sequence of {\it complete}
1-differentiable deformations of the Poisson bracket associated with the
symplectic form. 
This is what we prove in the present paper. In addition we are able to give
explicit formulas for these deformations.
As a corollary we show that the characteristic class of a star
product (in the formal Poisson bivector form, as in \cite{K} ) is explicitely
written in the expression of the star product.\\
\noindent Our ultimate goal is to show that one can reconstruct any star
product on a Poisson manifold from any other by adding (in a sense to define)
some 1-differentiable deformations. To find a precise definition of this
``addition" is related with the hope that, for any given star product on
a {\it Poisson ma\-ni\-fold}, the method used here can give a way to identify 
the formal Poisson bivector type characteristic class (see \cite{K}) of the star
product by looking solely at its explicit formula.

The paper is constructed as follows:
in Section 2 we define an algebraic operation (a contraction)
and give some useful formulas related to it. Section 3 is devoted to
the statement of our main result describing the impact of a
1-differentiable modification of a cochain at any given level $k$ on
subsequent levels, and developing the above mentioned consequences.
In Section 4 we give the proofs, relying for clarity of the
exposition on some intermediary lemmas which we also prove, omitting
details of straightforward computations. The last Section gives an idea about
the forms and the occurences of the 1-differentiable terms appearing in the 
formula of a Fedosov star-product. In the appendix we give the complete details
of the proofs.


\section{Definitions}

Throughout the article $M$ will be a symplectic manifold of dimension $2d$,
 $ \omega$ its non degenerate closed 2-form and
$\mu :\bigotimes^\bullet T^* M \longrightarrow \bigotimes^\bullet TM$
the canonical isomorphism given by $\omega$.
For $\alpha \in \Gamma (M,\otimes ^2 T^*M)$ we note
 $\bar{\alpha}= \mu(\alpha)$.\\
If $\alpha\in\Gamma(M,\bigwedge^2T^*M)\subset\Gamma(M,\otimes^2T^*M)$,
we write (in local coordinates on an arbitrary chart)
$\alpha=\alpha_{ij}dx^i\otimes dx^j=\frac{1}{2}\alpha_{ij}dx^i\wedge dx^j$
and the same for skewsymmetric bivectors. For the symplectic form
$\omega=\omega_{ij} dx^i\otimes dx^j = \frac{1}{2}\omega_{ij} dx^i\wedge dx^j$
we define $\bar\omega=\bar\omega^{ij} \partial_i\otimes \partial_j
= \frac{1}{2}\bar\omega^{ij} \partial_i\wedge \partial_j$ by
$\omega_{ij} \bar\omega^{jk} = \delta^k_j$.\\
 We have $\bar{\alpha}= \mu (\alpha) = -\bar{\omega}^{ir}
 \bar{\omega}^{js}\alpha _{rs}\partial _i \otimes \partial _j$ for $\alpha \in
\Gamma (M,\otimes ^2 T^*M)$.
 As in \cite{F1, F2} we use the Einstein convention on repeated 
indices $i,j = 1, \ldots , 2d$.

\begin{defi} The ``diamond" contraction.\\
$\diamond$ is the following contraction operation:
\begin{displaymath}
\begin{array}{llll}
    & \Gamma (M,\otimes ^2 T^*M)^{\otimes 2} & \longrightarrow & \Gamma
(M,\otimes ^2 T^*M)\\
    & \alpha \otimes \beta & \longmapsto & \alpha \diamond \beta
              = \bar{\omega}^{rs} \alpha_{ri} \beta_{sj} dx^i \otimes dx^j \\
    & & & \\
and & \Gamma (M,\otimes ^2 TM)^{\otimes 2} & \longrightarrow &
 \Gamma (M,\otimes ^2 TM)\\
    &  A \otimes B & \longmapsto & A \diamond B
                 = \omega_{rs} A^{ri} B^{sj} \partial_i \otimes \partial_j
\end{array}
\end{displaymath}
\noindent We define $\alpha^{\diamond n}=\alpha\diamond\alpha^{\diamond (n-1)}$
 and              $A^{\diamond n}=A\diamond A^{\diamond (n-1)}$.
\end{defi}

\medskip

\begin{prop}
Let $\alpha , \beta \in \Gamma (M,\bigwedge^2 T^*M)$ and
  $A, B \in \Gamma (M,\bigwedge^2 TM)$.\\
(i) The isomorphism $\mu$ acts multiplicatively with respect to the
diamond contraction:
$ \mu (\alpha \diamond \beta )= \mu (\alpha) \diamond \mu (\beta)$ ;
 $\mu^{-1}(A\diamond B) = \mu^{-1}(A) \diamond \mu^{-1}(B)$\\ 
(ii) In this case (values in skewsymmetric tensors) we have also\\
$\alpha^{\diamond n}\in \Gamma (M,\bigwedge ^2 T^*M)$ and 
 $A^{\diamond n}\in \Gamma (M,\bigwedge ^2 TM)$. \\ Moreover, for $l+m=n$, 
$ \alpha^{\diamond n}=\alpha^{\diamond l}\diamond\alpha^{\diamond m}$
and the same for $A$.

\end{prop}
Proof:\\ (i) direct computation.\\
         (ii) By induction. We use that, for 
$\alpha ,\beta ,\gamma \in \Gamma (M,\bigwedge ^2 T^*M)$, we have:\\
$\alpha\diamond (\beta\diamond\gamma)=(\beta\diamond\alpha)\diamond\gamma$ \ \
and for $\tau\in \Gamma (M,\otimes ^2 T^*M)$, \\  
$\omega\diamond\tau=-\tau\ \ ;(\tau\diamond\omega)_{ij}=\tau_{ji}$.  
The contravariant part is deduced by (i).
\begin{flushright} $\blacksquare $ \end{flushright}

One could also remark that, with the notations used in \cite{F1} \cite{F2}, 
in the Weyl bundle framework we have $\alpha^{\diamond 2}=
\frac{4i}{\hbar}(\delta^{-1}\alpha\circ\delta^{-1}\alpha)$.\\   
\noindent For $\alpha \in \Gamma (M,\bigwedge ^2 T^*M)$, we define
the coefficients $\alpha^{\diamond n}_{ij}$ by
$$\alpha^{\diamond n}=\alpha^{\diamond n}_{ij}dx^i \otimes dx^j
		=\frac{1}{2}\alpha^{\diamond n}_{ij}dx^i \wedge dx^j$$


\section{Results}

Let $(M,\omega)$ be a symplectic manifold, $\nabla$ a symplectic connexion 
on $M$ and $R$ its curvature. Let $*$ be a Fedosov star product on
$M$: for $f,g \in \mathcal C^\infty (M)$,
\begin{equation}
f*g=f.g-\frac{i\hbar}{2} \bar{\omega}(f,g) + \sum_{n\geq 1} \hbar^n C_n(f,g) \ .
\end{equation}

\noindent Let $\Omega \in \omega + \hbar Z_{DR}^{2}(M)[[\hbar]]$ be the
\textit{Weyl curvature} of $*$ (actually Fedosov \cite{F1,F2} takes $-\Omega$),
where $Z_{DR}^{2}(M)$ is the space of de Rham 2-cocycles on $M$. We denote by
$H_{DR}^{2}(M)$ the second de Rham cohomology space.
The \textit{characteristic class} of $*$ is the class of
$\Omega$ in\\ $[\omega] + \hbar H_{DR}^{2}(M)[[\hbar]]$.

Let $\tilde{*}$ be the Fedosov star product of Weyl curvature $\tilde{\Omega} =
 \Omega + \hbar^k \alpha$, with $\alpha \in  Z_{DR}^{2}(M)$
\begin{equation}
 f\tilde{*}g=f.g-\frac{i\hbar}{2} \bar{\omega}(f,g) +
\sum_{n\geq 1} \hbar^n \tilde{C}_n(f,g) \ . 
\end{equation}
We know by \cite{BCG} that $C_n = \tilde{C}_n $ for all $n\leq k$ and
$\tilde{C}_{k+1} = C_{k+1} + \frac{i}{2} \bar{\alpha}$.
What happens for $n>k+1$? This is the subject of the following proposition.

\begin{prop}  \label{p}
The change $\Omega \rightarrow \Omega + \hbar^k \alpha$ adds the series 
$\sum_{ p\geq 1} \hbar^{pk+1} \bar{\alpha}^{\diamond p}$ to the explicit
expression of the star-product $*$.
The series $\bar{\omega}-\sum_{ p\geq 1} \hbar^p \bar{\alpha}^{\diamond p}$
is a formal Poisson bracket and contains all the 1-differentiable terms of 
$\tilde{*}$ not depending explicitely on $R, \ \Omega-\omega$ and derivatives 
of $\alpha$.

\end{prop}

\smallskip

\noindent \textbf{Remarks:}
\begin{enumerate}
\item Let us choose $\Omega = \omega + \hbar^k \alpha$.
 We define the skewsymmetric bivector $\bar\Omega$ by 
$\Omega_{ij}\bar\Omega^{jk} = \delta_i^k$ and we obtain by easy computations
\begin{equation}
 \bar\Omega = 
\bar{\omega}-  \sum_{ p\geq 1} ( \hbar^{k} \bar{\alpha})^{\diamond p}=
\bar{\omega}-  \sum_{ p\geq 1}  \hbar^{pk} \bar{\alpha}^{\diamond p}
\end{equation}
 Since $d_{DR}\Omega =0$ we have $[\bar\Omega ,\bar\Omega ]=0$
for the Schouten bracket. Thus $\bar\Omega$ is a formal Poisson bivector.

\item Formula (3) is valid for a formal $\alpha$, i.e.
$ \bar\Omega =\bar{\omega}-  \sum_{ p\geq 1}
(\bar{\alpha}^\hbar )^{ \diamond p}$, 
with $\alpha^{\hbar}=\hbar^{k_1} \alpha_1 +\hbar^{k_2} \alpha_2
 + \cdots + \hbar^{k_m} \alpha_m + \cdots, \ \ \Omega = \omega +  \alpha^\hbar$.
This gives all the 1-differentiable deformations of $\bar\omega$.
If $[\alpha^\hbar]=[\beta^\hbar]$ in $H^2(M)[[\hbar]]$ the two resulting
formal Poisson brackets are equivalent (see \cite{K}).

\end{enumerate}
\medskip

Proposition \ref{p} can be reformulated in order to make more explicit 
the relation between this formal Poisson bracket and the characteristic class
of the star product:

\begin{cor} 
Let $*$ be the Fedosov star product of trivial Weyl curvature $\omega$ 
and $\tilde{*}$ the one with curvature $\Omega = \omega +  \alpha^\hbar$.
$\Omega$ appears explicitely in the formal Poisson form $\bar{\Omega}$ 
as a part of the formula for $\tilde{*}$. $\bar{\Omega}$ can be seen
as all the 1-differentiable terms of $\tilde{*}$ not containing $R$ or
derivatives of $\alpha^\hbar$. We have:
\begin{eqnarray}
f\tilde{*}g&=&f*g+\frac{i\hbar}{2}\sum_{ p\geq 1}
(\bar{\alpha}^\hbar )^{ \diamond p}(f,g)+ \rho (f,g) \nonumber \\
  &=& f.g - \frac{i\hbar}{2}\left(\bar\omega (f,g)
          -\sum_{ p\geq 1} (\bar{\alpha}^\hbar )^{ \diamond p}(f,g)\right)
              +\sum_{n\geq 2} \hbar^n C_n(f,g)+ \rho (f,g)\nonumber \\
&=& f.g -\frac{i\hbar}{2}\bar\Omega (f,g)+\sum_{n\geq 2} \hbar^n C_n(f,g)
                                                                  +\rho (f,g)
\end{eqnarray}
where the terms occuring in the remainder $\rho$ either depend explicitely 
on the curvarture $R$ or on derivatives of $\alpha^\hbar$, 
or are not 1-differentiable.

\end{cor}

\medskip

\noindent \textbf{Remark:} Expression (4) shows that on $(M,\omega)$,
for a star-product $\tilde{*}$ of characteristic class $\Omega$, the
corresponding bracket
$\{f,g\}_{\tilde{*}}=\frac{i}{\hbar}(f\tilde{*}g - g\tilde{*}f)$ can be seen 
not only as a deformation of the Lie algebra $(\mathcal{C}(M),\bar\omega )$ 
but also as	a deformation of the ``formal" Lie algebra 
$(\mathcal{C}(M)[[\hbar]],\bar\Omega )$.
Thus $\tilde{*}$ can be viewed as the star-product of trivial
characteristic class on the ``formal" symplectic manifold $(M, \Omega)$, i.e.
$M$ endowed with the formal symplectic structure given by $\Omega$.
We can also consider $\tilde{*}$ as a deformation of $*$ with 
$\frac{i}{2}\sum_{ p\geq 1}(\bar{\alpha}^\hbar )^{ \diamond p}$ as 
infinitesimal deformation.


\section{Proofs}

\subsection {Fedosov notations}

We use the notations of Fedosov \cite{F1,F2}:
we choose a symplectic connexion $\nabla$ on the symplectic manifold
$M$ and we denote by $\partial$ the covariant exterior derivative associated 
to $\nabla$. Let $\circ$\ ,\ $[.\ ,\ .]$ be, respectively, the (Moyal) product 
and the bracket on ${\cal W}=\Gamma (M,W)$ , the sections
of the Weyl bundle $W$ associated to $M$, and 
$\delta$,$\delta^{-1}$ the operators on ${\cal W}$ defined in \cite{F1,F2}.
We construct on $W$ an Abelian connexion
$D=\partial -\delta +\frac{i}{\hbar} [r,.]$. $\Omega$ is the curvature of $D$.
Defining $Q=R+(\Omega -\omega)$, $r$ is the unique solution, under suitable 
conditions, of the equation:
$$r=\delta^{-1}Q + \delta^{-1}(\partial r +\frac{i}{\hbar}r^2)$$
Fedosov shows that the space 
${\cal W}_D=\{a=a(x,y,\hbar)\in {\cal W}\ |\ Da=0\}$
is isomorphic to
${\cal{C}}^\infty (M)[[\hbar]]$ by the iso\-mor\-phism $\sigma_D$ 
given by the equation
\begin{equation}
\sigma_D^{-1}(f) = a = f+\delta^{-1}(\partial a+\frac{i}{\hbar}[r,a])\ ,
 \ f\in \mathcal C^\infty (M)[[\hbar]] \ .
\end{equation}
$\sigma_D$ is the restriction to ${\cal W}_D$ of the projection $\sigma$,
$\sigma (a) = a(x,0,\hbar)$ ($\sigma$ 
replaces the $y$'s by 0).
So the star product corresponding to $\Omega$ is given by 
$f*g=\sigma_D (\sigma_D^{-1}(f)\circ\sigma_D^{-1}(g))$.
For $a \in  {\cal W}$,  $a^{(n)}$ will denote the part of degree $n$ of $a$ 
in the usual filtration of the Weyl bundle while $a_n$ will be defined as $a$ 
modulo the terms of degree $>n$, i.e. $a_n =a^{(0)}+a^{(1)}+\cdots+a^{(n)}$.

\subsection{Equation for $\tilde{r}$}

We are now looking for an $\tilde{r}$ giving an Abelian connexion $\tilde{D}$
of Weyl curvature 
$\tilde{\Omega}=\Omega + \hbar^k\alpha=\frac{1}{2}\Omega_{ij} dx^i\wedge dx^j
+\hbar^k\frac{1}{2}\alpha_{ij} dx^i\wedge dx^j$.
So $\tilde{r}$ satisfies \\
$\tilde{r}=\delta^{-1}Q +\delta^{-1}(\hbar^k\alpha)
+ \delta^{-1}(\partial \tilde{r} +\frac{i}{\hbar}\tilde{r}^2)$.
This modification gives some additional terms of interest:

\begin{itemize}

\item the first appears in degree $2k+1$. We have
$\tilde{r}^{(2k+1)}=r^{(2k+1)} + \hbar^k s_1$ with
$s_1 = \delta^{-1} \alpha =\frac{1}{2}\alpha_{ij}y^i dx^j
=\sigma_1\ \alpha_{ij}y^i dx^j$

\item in degree $4k+1$ :
$$\hbar^{2k}s_2 = \hbar^{2k} \delta^{-1}(\frac{i}{\hbar}(\delta^{-1}\alpha
\circ\delta^{-1}\alpha))= \hbar^{2k} \delta^{-1}(\frac{i}{\hbar}(s_1\circ s_1))
$$ $$=\hbar^{2k} \frac{1}{8} \alpha^{\diamond 2}_{ij}y^i dx^j=
\hbar^{2k}\sigma_2\ \alpha^{\diamond 2}_{ij}y^i dx^j$$

\end{itemize}
Actually we consider only the terms depending exclusively on the $s_i$'s.
So the next degrees to consider are (using $2\tilde{r}=[\tilde{r}, \tilde{r}]$):
\begin{itemize}
\item degree $6k+1$ :
$$\hbar^{3k}s_3=\hbar^{3k} \delta^{-1}\frac{i}{2\hbar}([s_1, s_2]+[s_2, s_1])
=\hbar^{3k} \sigma_3\ \alpha^{\diamond 3}_{ij}y^i dx^j$$

\item . . .

\item degree $2pk+1$ : $$\hbar^{pk}s_p=\hbar^{pk} \delta^{-1}\frac{i}{2\hbar}
\sum_{l+m=p;\ l,m\geq1}[s_l, s_m]
=\hbar^{pk} \sigma_p\ \alpha^{\diamond p}_{ij}y^i dx^j.$$
\end{itemize}

\noindent The $\sigma_p$'s will be described in the proof of Lemma \ref{coef}.

\noindent In these computations we have used the 
fol\-lo\-wing straight\-for\-ward lem\-ma:

\begin{lemma} 
\ \ $\delta^{-1}\left(\frac{i}{2\hbar}[\alpha^{\diamond l}_{ij}y^i dx^j,
\alpha^{\diamond m}_{ij}y^idx^j]\right)
=\frac{1}{2}\alpha^{\diamond (l+m)}_{ij}y^i dx^j $.
\end{lemma}

\subsection{Equation for $\tilde{a}$}

We now describe the consequences of these changes in the computations of
$\sigma_{\tilde{D}}^{-1}(f) = \tilde{a}
=f+\delta^{-1}(\partial\tilde{a}+
\frac{i}{\hbar}[\tilde{r},\tilde{a}])$, \ $f\in \mathcal C^\infty (M)$.
Recall that for any Weyl curvature, $a^{(1)}=\partial_jf\ y^j.$
We then obtain:
\begin{itemize}

\item in degree $2k+1$:
\begin{eqnarray*}
\tilde{a}^{(2k+1)}&=&
a^{(2k+1)}+\hbar^k\delta^{-1}(\frac{i}{\hbar}[s_1,a^{(1)}]=
a^{(2k+1)}+\hbar^k x_1 \\
&=&a^{(2k+1)} +
\frac{1}{2} \hbar^k \bar\omega^{il}\alpha_{ij}\partial_lf\ y^j
\end{eqnarray*}

\item among other additional terms, there is (with $\varkappa_2=\frac{3}{2}$):
\begin{eqnarray*}\hbar^{2k}x_2&=&
\hbar^{2k}\delta^{-1}(\frac{i}{\hbar}([s_2,a^{(1)}] + [s_1,x_1])\\
&=&\hbar^{2k}
\varkappa_2 \bar\omega^{il}(\alpha^{\diamond 2})_{ij}\partial_lf\ y^j. 
\end{eqnarray*}

\item . . .

\item in degree $2pk+1$\ (with $x_0=a^{(1)}$ , $\varkappa_0=1$): 
\begin{eqnarray*} \hbar^{pk}x_p&=&
\hbar^{pk}\delta^{-1}(\frac{i}{\hbar}\sum_{l+m=p;\ l\geq1,m\geq0}[s_l,x_m])\\
&=&\hbar^{pk}\varkappa_p\
\bar\omega^{il}(\alpha^{\diamond p})_{ij}\partial_lf\ y^j.
\end{eqnarray*}

\end{itemize}

\noindent The $\varkappa_p$'s will be described in the proof of 
Lemma \ref{coef}.

\noindent In these computations we have used the following lemma:

\begin{lemma} \ \
$\delta^{-1}(\frac{i}{\hbar}[\alpha^{\diamond l}_{ij}y^i dx^j\
,\ \bar\omega^{il}\alpha^{\diamond m}_{ij}\partial_lf\ y^j])
= \bar\omega^{il}\alpha^{\diamond (l+m)}_{ij}\partial_lf\ y^j$.
\end{lemma}

\subsection{End of the proof of Proposition \ref{p} }

Finally, from the formula $$\tilde{C}^n(f,g)
=\sigma\left( (\tilde{a}\circ\tilde{b})^{(2n)}\right)
 = \sigma\left(\sum_{l+m=2n}\tilde{a}^{(l)}\circ\tilde{b}^{(m)}\right) ,$$ with
 $\tilde{a}=\sigma_{\tilde{D}}^{-1}(f),\
\tilde{b}=\sigma_{\tilde{D}}^{-1}(g),\ f,g \in \mathcal C^\infty (M)$ and
taking into account (straightforward computations)

\begin{lemma} 
$$\sigma\left(x^{(a)}_l\circ x^{(b)}_m \right)
= \frac{i \hbar}{2}\varkappa_l\varkappa_m (\bar{\alpha}^{\diamond(m+l)})^{ij} 
\partial_if\partial_jg
= \frac{i \hbar}{2}\varkappa_l\varkappa_m \bar{\alpha}^{\diamond(m+l)}(f,g)$$
\end{lemma}
 we obtain:

\begin{itemize}
\item $\tilde{C}^{k+1}(f,g)=C^{k+1}(f,g) + \frac{i}{2}\ \bar\alpha (f,g)$
\item $\tilde{C}^{2k+1}(f,g)=C^{2k+1}(f,g) + i\ c_2\
 \bar\alpha^{\diamond 2} (f,g) + \rho_2(f,g)$
\item . . .
\item $\tilde{C}^{pk+1}(f,g)=C^{pk+1}(f,g) + i\ c_p\
 \bar\alpha^{\diamond p} (f,g) + \rho_p(f,g)$
\end{itemize}


\begin{lemma} \label{coef} \ \ \
$c_p = \frac{1}{2}$
\end{lemma}

\textit{Proof of Lemma 4:} \\ Define $S(x)=\sum_{n\geq1}\sigma_nx^n$ and
$X(x)=\sum_{n\geq0} \varkappa_n x^n$.\\ Since
${\displaystyle \sigma_n=\frac{1}{2}\sum_{l+m=n;\ l,m\geq1}\sigma_l\sigma_m}$,
we find $\frac{1}{2}S^2(x)=S(x)-\frac{1}{2}x$ and therefore
$S(x)=-\sqrt{1 - x}+1$. \\
In the same way, since 
${\displaystyle \varkappa_n=\sum_{l+m=n;\ l\geq0,m\geq1} 
\varkappa_l \sigma_m}$, 
we have \\ $S(x) X(x) = X(x) - 1$ and so $X(x)=\frac{1}{\sqrt{1-x}}$.\\
Finally, since $ {\displaystyle c_n = \frac{1}{2} \sum_{l+m=n;
\ l,m\geq 0} \varkappa_l \varkappa_m\ , c_n}$ is the $n^{th}$ coefficient
of the Taylor expansion of $\frac{1}{2}X^2(x)= \frac{1}{2(1-x)}$.\\
This completes the proof of Lemma \ref{coef}. \hfill $\blacksquare $

\medskip

\noindent {\it Specificity of these terms:}

In Equation (5) for $r$, we considered all the terms
involving solely $\delta^{-1}\alpha$. The other terms always depend at least
on $\delta^{-1}Q = \delta^{-1}(R + \Omega - \omega)$ 
or $ (\delta^{-1}\partial)^n\delta^{-1}\alpha$, $n\geq1$, and therefore
involve $Q$ or derivatives of $\alpha$.\\
At the next step, the only 1-differentiable term of $a=\sigma_D^{-1}(f)$
constructed without $\delta^{-1}Q$ is $a^{(1)}= y^i \partial_if$.
So we have considered all the terms obtained inductively with 
$\frac{i}{\hbar}\delta^{-1}[\tilde{r},\tilde{a}]$ 
mixing the ones found in the first step and $a^{(1)}$. 
The terms coming from the part ``$\delta^{-1}\partial \tilde{a}$"
of the equation on $\tilde{a} $ will depend on derivatives of $\alpha$ or 
won't be 1-differentiable anymore.\\
The last step just contracts these selected terms with $\bar{\omega}$.
 \hfill $\blacksquare $ (Prop. \ref{p})

\medskip

Some ideas about the form and the propagation of the other 
1-dif\-fe\-ren\-tiable 
terms are given in the next section.


\subsection{Proof of the Corollary}

The corollary is straightforward for $\Omega=\omega+\hbar^k\alpha$,
$\alpha \in \Gamma (M,\bigwedge^2 T^*M)$.
For $\alpha^{\hbar}\in \Gamma (M,\bigwedge^2 T^*M)[[\hbar]]$
the proof of Proposition\ref{p} is easily adaptable. Indeed, let us take
$\alpha^{\hbar}=\hbar^{k_1} \alpha_1 +\hbar^{k_2} \alpha_2
 + \cdots + \hbar^{k_n} \alpha_n$
with $\alpha_q=\frac{1}{2}(\alpha_q)_{ij}dx^i\wedge dx^j$.
We are looking for a $\tilde{r}$ such that

\begin{equation} \label{r}
\tilde{r}=\delta^{-1}R+\delta^{-1}\alpha^{\hbar}
+\delta^{-1}(\partial\tilde{r}+\frac{i}{\hbar}\tilde{r}^2)
\end{equation}
So $\tilde{r}_{2k_n+1}=r_{2k_n+1}+\delta^{-1}\alpha^{\hbar}+\cdots $.
We can also write $\tilde{r}=r+\delta^{-1}\alpha^{\hbar}+\cdots $\\ so that
another application of (\ref{r}) gives:
$$\tilde{r}=\delta^{-1}R+\delta^{-1}\alpha^{\hbar}+\delta^{-1}
\left(\partial\tilde{r}+\frac{i}{\hbar}(r+\delta^{-1}\alpha^{\hbar}
+\cdots)^2\right),$$
and we see that $\tilde{r}$ contains
$\delta^{-1}(\frac{i}{\hbar}(\delta^{-1}\alpha^{\hbar}\circ\delta^{-1}
\alpha^{\hbar}))=\frac{1}{8}(\alpha^{\hbar})^{\diamond 2}_{ij}y^idx^j$
 (actually this is true from $\tilde{r}_{4k_n+1}$).
Iteration of this process gives way to the same computations as before.\\
The same argument can be used for solving the equation
\begin{equation}
\tilde{a}=f+\delta^{-1}(\partial\tilde{a}+\frac{i}{\hbar}
[\tilde{r},\tilde{a}]).
\end{equation}
 We find $\tilde{a}=
a+\sum_{p\geq 1}\varkappa_p\widehat{(\alpha^{\hbar}_f)^{\diamond p}}+\cdots $
with $\widehat{(\alpha^{\hbar}_f)^{\diamond p}}
=\omega^{il}((\alpha^{\hbar})^{\diamond p})_{ij}\partial_lf\ y^j$. 

\smallskip

\noindent Then
\begin{eqnarray}
f\tilde{*}g&=&\sigma(\tilde{a}\circ\tilde{b})\nonumber \\
&=&\sigma (a\circ b)+\sigma{\left(\sum_{l\geq 1}\varkappa_l
\widehat{(\alpha^{\hbar}_f)^{\diamond l}}\
\circ \ \sum_{m\geq 1}\varkappa_m\widehat{(\alpha^{\hbar}_g)^{\diamond m}}
\right)}+\rho(f,g) \nonumber \\
&=&f*g+\frac{i\hbar}{2}\sum_{ p\geq 1} (\bar{\alpha}^\hbar )^{ \diamond p}
(f,g)+\rho (f,g)
\end{eqnarray}
by lemmas 3 and 4.\\
So the corollary is proved for 
$\alpha^{\hbar} \in \Gamma (M,\bigwedge^2 T^*M)[\hbar]$ and
by induction for $\alpha^{\hbar} \in \Gamma (M,\bigwedge^2 T^*M)[[\hbar]] $.


\section{Ideas about the form and the propagation of 1-differentiable terms}

In this section we want to give an idea about the occurences and the forms
of the other 1-differentiable terms that can appear in the explicit expression 
of a Fedosov star-product.
We use the notations of Section~4.

\subsection{In the star-product of Weyl curvature $\Omega=\omega$ }

$\delta^{-1}R\circ\delta^{-1}R$ contains a term without any $"y"$, so it is
a 2-form on $M.$ Let's denote it $\beta_0$. $\hbar^{2}\delta^{-1}\beta_0$
appears in $r^{(5)}$. More generally let's denote by
$\beta_n$ the 2-form (part without any $"y"$) appearing in
$ \left( (\delta^{-1}\partial)^n\delta^{-1}R\right)^2$, i.e. we have
$\frac{i}{\hbar}
\left( (\delta^{-1}\partial)^n\delta^{-1}R\right)^2|_{y=0} 
=\hbar^{n+2}\beta_n$ which appears in $r^{(2n+5)}$.
For $n$ odd, $\beta_n = 0$. Then 
\begin{equation}
\sigma\left(\delta^{-1}\frac{i}{\hbar}[\hbar^{n+2}\delta^{-1}\beta_n,a^{(1)}]
\circ b^{(1)} + a^{(1)}\circ \delta^{-1}\frac{i}{\hbar}
[\hbar^{n+2}\delta^{-1}\beta_n,b^{(1)}]\right)
\end{equation}
gives a 1-differentiable term in every $C_{3+n}\ , \ n \in 2\mathbb{N}$.
Since, between the part ``in $a$" and the one ``in $b$", it uses the product
$\circ$ only at the first order in $\hbar$, it is skewsymmetric.
\begin{equation}
\sigma\left(\delta^{-1}\frac{i}{\hbar}[(\delta^{-1}\partial)^n\delta^{-1}R,
a^{(1)}]\circ \delta^{-1}\frac{i}{\hbar}[(\delta^{-1}\partial)^n\delta^{-1}R,
b^{(1)}]\right)
\end{equation}
gives a symmetric (resp. skewsymmetric) 1-differentiable term in $C_{3+n}$
for $n$ odd (resp. for $n$ even).

\smallskip

So, at worse, each $C_l$ contains a 1-differentiable part for $l\geq3$. \\
But the above term (10) might be cancelled because it might appear under 
other forms. For example, in the case $\beta_0$,
three kinds of 1-differentiable terms appear in $C_3$:
\begin{itemize}
\item[1)] $\sigma\left(\delta^{-1}\frac{i}{\hbar}[\delta^{-1}R,a^{(1)}]\circ 
\delta^{-1}\frac{i}{\hbar}[\delta^{-1}R,b^{(1)}]\right)$

\item[2)] $\begin{array}{ll} \sigma\Big(\delta^{-1}\frac{i}{\hbar}[\delta^{-1}R&,
\delta^{-1}\frac{i}{\hbar}[\delta^{-1}R,a^{(1)}]]\circ b^{(1)}\\
&+a^{(1)}\circ \delta^{-1}\frac{i}{\hbar}[\delta^{-1}R,
\delta^{-1}\frac{i}{\hbar}[\delta^{-1}R,b^{(1)}]]\Big)
\end{array}$

\item[3)] $\sigma\left(\delta^{-1}\frac{i}{\hbar}[\hbar^{2}\delta^{-1}\beta_0,
a^{(1)}]\circ b^{(1)} + a^{(1)}\circ \delta^{-1}\frac{i}{\hbar}
[\hbar^{2}\delta^{-1}\beta_0,b^{(1)}]\right)$

\end{itemize}

In this case these three terms are the same, up to a positive coefficient,
so they cannot cancel. I do not know if this always happens.
{\it This kind of phenomena can occur for all the terms we consider 
in these sections}.

\medskip

One can observe an interesting phenomenon of propagation: \\ since 
$\hbar^{n+2}\delta^{-1}\beta_n$ appears in $r$, it propagates
exactly in the same way as $\hbar^k\delta^{-1}\alpha$ in the
proof of Proposition \ref{p}, so a series 
${\displaystyle \sum_{ p\geq 1} \hbar^{p(n+2)+1} 
\bar{\beta_n}^{\diamond p}(df,dg)}$
appears.

\subsection{Effects of the change $\Omega \rightarrow \Omega+\hbar^k\alpha$}

\subsubsection{Mixed terms} (by ``mixed terms" we mean the terms 
involving both $R$ and $\alpha$).

\smallskip

It is not difficult to see that in $\tilde{C}_{k+2}$ there are no
supplementary 1-differentiable terms compared to $C_{k+2} \ (k\geq 2)$.\\
In $\tilde{C}_{k+3}$, there is one. It is the case $n=0$ of the
following fact: the term
\begin{displaymath} 
\begin{array}{ll}
\sigma\Big(\delta^{-1}\frac{i}{\hbar}[\delta^{-1}\hbar^k\alpha,&
\delta^{-1}\frac{i}{\hbar}
[\delta^{-1}\hbar^{n+2}\beta_n,a^{(1)}]]\circ b^{(1)} \\
 &+a^{(1)}\circ \delta^{-1}\frac{i}{\hbar}[\delta^{-1}\hbar^k\alpha,
\delta^{-1}\frac{i}{\hbar}[\delta^{-1}\hbar^{n+2}\beta_n,b^{(1)}]]\Big)
\end{array} 
\end{displaymath}
is 1-differentiable, skewsymmetric and part of $\tilde{C}_{k+3+n}\ ,
\ \forall n$ even.\\
And
\begin{displaymath} 
\begin{array}{ll}
\sigma\Big(\delta^{-1}\frac{i}{\hbar}[\delta^{-1}\hbar^k\alpha,&
\delta^{-1}\frac{i}{\hbar}[(\delta^{-1}\partial)^n\delta^{-1}R, a^{(1)}]]\\ 
&\circ \ \delta^{-1}\frac{i}{\hbar}[\delta^{-1}\hbar^k\alpha,
\delta^{-1}\frac{i}{\hbar}[(\delta^{-1}\partial)^n\delta^{-1}R,b^{(1)}]]\Big)
\end{array} 
\end{displaymath}
is 1-differentiable, symmetric for $n$ odd, skewsymmetric for $n$ even.\\
So the change $\Omega \rightarrow \Omega + \hbar^k\alpha$
can give a supplementary
1-differentiable term in every $\tilde{C}_l \ f\!or\ l\geq k+3$.

\medskip

\noindent Another phenomenon of propagation can be observed: \\ denote 
$X_n(f) = \frac{i}{\hbar}[\delta^{-1}\beta_n,a^{(1)}] 
= \bar\omega^{il}\beta_{n,ij}\partial_lf dx^j$. It is a 1-form on $M$.
In the proof of Proposition \ref{p}, it is possible to replace 
$a^{(1)}=\delta^{-1}\partial f$ by $\delta^{-1}X_n(f)$. Thus the 
characteristic class can appear again in the form 
$\sum_{ p\geq 1} \hbar^{kp+2n+5} \bar{\alpha}^{\diamond p}(X_n(f),X_n(g))$.

\subsubsection{Terms purely in $\alpha$}

The first 1-differentiable terms not involving $R$ and
depending on deri\-va\-ti\-ves of $\alpha$ can appear in $\tilde{C}_{2k+2}$:
in the same way than in the preceding subsection, 
for $n\geq1$,
$$ \sigma\left(\delta^{-1}\frac{i}{\hbar}
[\hbar^k(\delta^{-1}\partial)^n\delta^{-1}\alpha,a^{(1)}]
\circ 
\delta^{-1}\frac{i}{\hbar}
[\hbar^k(\delta^{-1}\partial)^n\delta^{-1}\alpha,b^{(1)}]\right)$$
is 1-differentiable, symmetric for $n$ odd, skewsymmetric for $n$ even.

\smallskip

\noindent $\frac{i}{\hbar}\left( (\delta^{-1}\partial)^n
\delta^{-1}\alpha\right)^2 |_{y=0} = 
\hbar^{2k+n}\gamma_n$ where $\gamma_n$ is a 2-form on $M$. $\gamma_n = 0$
for $n$ odd. As with $\beta_n$ (eq. 9) one can construct a skewsymmetric 
1-diffe\-ren\-tia\-ble term with $\gamma_n$.

\smallskip

\noindent So one can find 1-differentiable terms of these types in every
$\tilde{C}_{2k+1+n}$, $n\geq1$.

\bigskip

Defining
$Y_n(f) = \frac{i}{\hbar}[\delta^{-1}\hbar^{2k+n}\gamma_n,a^{(1)}]$,
$\sum_{ p\geq 1} \hbar^{p(2k+n)+1} \bar{\gamma}_n^{\diamond p}(df,dg)$ and
$\sum_{ p\geq 1} \hbar^{kp+4k+2n+1} \bar{\alpha}^{\diamond p}(Y_nf,Y_ng)$
appear in the formula.

\bigskip

\noindent \textbf{Remarks:}
\begin{enumerate}
\item For simplicity in 5.2.1 we have considered 1-differentiable 
``mixed" terms not involving derivatives of $\alpha$, but there exist, 
for example, terms like $\beta_n^{\diamond p}(Y_m(f),Y_m(g))$.

\item It is also possible to have an idea of the propagation of
all the terms in $\alpha$, {\it not necessarily 1-differentiable},
which do not contain $R$. Let $u$ be the solution of 
$u=\hbar^k\delta^{-1}\alpha + \delta^{-1}(\partial u +
\frac{i}{\hbar}u^2)$ in $\Gamma (M,W\otimes T^*M)$ 
and $a_u$ the solution of $a_u=f+\delta^{-1}(\partial a_u +
\frac{i}{\hbar}[u,a_u])$ in $ \Gamma (M,W)$.
These solutions exist and are unique because 
$\delta^{-1}(\partial . \ +\frac{i}{\hbar}\ .\circ .)$ and 
$\delta^{-1}(\partial . \  + \frac{i}{\hbar}[u,.])$ raise degree 
(see \cite{F2}).
Putting $R=0$ in the expression of $\tilde{*}$ there is only
$\sigma(a_u\circ b_u)$ left. Formally, this is the expression of 
the Fedosov star-product on $\mathbb{R}^{2n}$ of Weyl curvature 
$\omega + \hbar^k \alpha$. 

For a Poisson bivector field $\pi$ we denote 
$$\pi^n(f,g) = \pi^{i_1j_1}\pi^{i_2j_2}\cdots\pi^{i_nj_n}
(\partial_{i_1}\cdots\partial_{i_n}f) (\partial_{i_1}\cdots\partial_{i_n}g)$$
So $\sigma(a_u\circ b_u)$ contains
$\exp(\frac{-i\hbar}{2}\bar\omega)(f,g)$ \cite{X}
and the other terms form the part ``purely" in $\alpha$ that is added to the
formula of the star-product when we change $\Omega=\omega$ in 
$\Omega=\omega+\hbar^k\alpha$. The conjecture is that this part contains 
$\exp(\frac{-i\hbar}{2}\bar\Omega)(f,g)$ with $\Omega=\omega+\hbar^k\alpha$.
Actually Proposition \ref{p} shows that it is true at order 1 
of differentiation.

\end{enumerate}

\bigskip

\centerline{\bf Acknowledgements}

\smallskip

I want to thank M. Flato, P. Gautheron and D. Sternheimer for asking me 
questions which push me to do this work (and especially D.S. for constant 
disponibility throughout the elaboration). I also want to thank F. Bidegain,
P. Bieliavsky and participants at the Warwick symposium in December 97 for 
numerous comments.

\bigskip

\bigskip

\pagebreak

\appendix

\centerline{\large\bf APPENDIX}

\bigskip

We give here the complete details of the proofs of the above results. After
this paper was completed, we received \cite{GR} and noticed (S. Gutt, private
communication) that, in a nonexplicit form, results similar to ours can be 
derived from there.

\section{Definitions}

Let $W$ be the ``Weyl" bundle i.e. the bundle of formal 
Weyl algebras defined in
\cite{F1,F2}. A section $a$ of $W$ is a sum of ``monomials" of the form
$\hbar^k a_{k;i_1,\ldots ,i_p}(x,y,\hbar)y^{i_1} \ldots y^{i_p} $. We 
give to it the degree $2k+p$ and this gives a filtration on 
${\cal{W}}=\Gamma(M,W)$.
For $a, b \in {\cal{W}}$  we have the following product:

\begin{eqnarray*}
a\circ b &=& \sum_{k=0}^\infty \left(-\frac{i\hbar}{2}\right)^k\frac{1}{k!}\ 
\bar\omega^{i_1 j_1} \ldots \bar\omega^{i_k j_k}\frac{\partial^k a}
{\partial y^{i_1}\ldots \partial y^{i_k}}  \frac{\partial^k b}
{\partial y^{j_1}\ldots \partial y^{j_k}}\\
&=& \sum_{k=0}^\infty a\circ_k b
\end{eqnarray*}

This product can be extended to the differential forms with values 
in $W$ by means of
the exterior product on the ``$dx^i$'s".

A graded commutator is defined by $[a,b]=a\circ b - (-1)^{q_1q_2} b\circ a$\ ,
\ for $ a\in \Gamma(M,W\otimes \bigwedge^{q_1}T^*M) $ and
$ b\in \Gamma(M,W\otimes \bigwedge^{q_2}T^*M) $.

We use the following two operators on the forms:
$$\delta a = dx^k\wedge \frac{\partial a}{\partial y^k} \ \ , \ \ 
\delta^{-1} a = \frac{1}{p+q} y^k i(\frac{\partial }{\partial x^k})a$$
for $a \in \Gamma(M,W\otimes \bigwedge^pT^*M)$ and of degree $q$ in
the filtration of $\cal{W}$.

\begin{lem}
$a,b \in \Gamma(M,W),\ \alpha,\beta \in \Gamma(M,W\otimes T^*M)$
$$(i) \ a\circ_k b = (-1)^k b\circ_k a \ ; 
\ \alpha\circ_k b = (-1)^k b\circ_k \alpha\ ;
\ \alpha\circ_k \beta  = (-1)^{k+1}\beta  \circ_k \alpha$$
$$(ii) [a,b]=2\sum_{p\geq 0}a \circ_{2p+1} b \ ;
\ [\alpha,b]=2\sum_{p\geq 0}\alpha \circ_{2p+1} b \ ;
\ [\alpha,\beta]=2\sum_{p\geq 0}\alpha \circ_{2p+1} \beta $$
\noindent In particular, $[a,b]=2a\circ_1b\ , \ [\alpha,b]=2\alpha\circ_1b\ ,\
[\alpha,\beta]=2\alpha\circ_1\beta$ for $\alpha, \beta, b$ of degree 1 in $y$. 
\end{lem}

\noindent {\bf Proof:}
\begin{eqnarray*}
\hbox{(i)} \   \ a\circ_k b &=& \left(-\frac{i\hbar}{2}\right)^k\frac{1}{k!}\ 
\bar\omega^{i_1 j_1} \ldots \bar\omega^{i_k j_k}\frac{\partial^k a}
{\partial y^{i_1}\ldots \partial y^{i_k}}  \frac{\partial^k b}
{\partial y^{j_1}\ldots \partial y^{j_k}}\\
&=& \left(-\frac{i\hbar}{2}\right)^k\frac{1}{k!}\ 
\bar\omega^{j_1 i_1} \ldots \bar\omega^{j_k i_k}\frac{\partial^k a}
{\partial y^{j_1}\ldots \partial y^{j_k}}  \frac{\partial^k b}
{\partial y^{i_1}\ldots \partial y^{i_k}}\\
&=& \left(-\frac{i\hbar}{2}\right)^k\frac{1}{k!}\ (-1)^k 
\bar\omega^{i_1 j_1} \ldots \bar\omega^{i_k j_k}
\frac{\partial^k b}{\partial y^{i_1}\ldots \partial y^{i_k}}
\frac{\partial^k a}{\partial y^{j_1}\ldots \partial y^{j_k}} \\ 
& & \hbox{as}\ \bar\omega^{ij}=-\bar\omega^{ji}\\
&=& (-1)^k \ b\circ_k a
\end{eqnarray*}

\medskip

$\alpha,\beta \in \Gamma(M,W\otimes T^*M)$ so we can write 
$\alpha=\alpha_i(x,y,\hbar) dx^i$ and $\beta = \beta_j(x,y,\hbar) dx^j$. 
$$\hbox{Then} \ \ \ \
\alpha\circ_k b=\alpha_i\circ_k b\ dx^i = (-1)^k b\circ_k \alpha_i dx^i
=(-1)^k b\circ_k \alpha$$ 
\begin{eqnarray*}
\hbox{and} \ \ \alpha\circ_k \beta &=&\alpha_i\circ_k \beta_j \ dx^i\wedge dx^j 
= (-1)^k\ \beta_j \circ_k \alpha_i \ dx^i\wedge dx^j\\ &=& 
(-1)^{k+1} \beta_j \circ_k \alpha_i\ dx^j\wedge dx^i
=(-1)^{k+1} b\circ_k \alpha
\end{eqnarray*}

\begin{eqnarray*}
\hbox{(ii)} \ \ [a,b]&=&a\circ b - b\circ a 
= \sum_{k\geq0}(a\circ_k b - b\circ_k a)
=\sum_{k\geq0} (1-(-1)^k)\ a\circ_k b \\ &=& 2 \sum_{p\geq0} a\circ_{2p+1} b
\end{eqnarray*}
\noindent $[\alpha,b]=\alpha\circ b - b\circ \alpha$ 
and $[\alpha,\beta]=\alpha\circ \beta + \beta \circ \alpha$ 
and the computations are the same.

\hfill $\blacksquare $

\bigskip

\noindent {\bf About the diamond product:}

\bigskip

Let $\alpha , \beta ,\gamma \in \Gamma (M,\bigwedge^2 T^*M)$,
$\tau\in \Gamma (M,\otimes ^2 T^*M)$.
We will show
\begin{enumerate}
\item
$ \mu (\alpha \diamond \beta )= \mu (\alpha) \diamond \mu (\beta)$

\item
$\alpha^{\diamond n}\in \Gamma (M,\bigwedge ^2 T^*M)$ and for $l+m=n$, 
$ \alpha^{\diamond n}=\alpha^{\diamond l}\diamond\alpha^{\diamond m}$
with the help of
\begin{enumerate}
\item
$(\alpha\diamond\beta)\diamond\gamma=\beta\diamond(\alpha\diamond\gamma)$ 

\item
$(\tau\diamond\omega)_{ij}=\tau_{ji}$

\end{enumerate}
\end{enumerate}

\noindent {\bf Proofs:}

\medskip

\noindent 1. 
\begin{eqnarray*}
(\mu (\alpha) \diamond \mu (\beta))^{ij} &=& 
\omega_{rs} \mu (\alpha)^{ri} \mu (\beta)^{sj}\\
&=& \omega_{rs}(-\bar{\omega}^{kr} \bar{\omega}^{li}\alpha _{kl})
(-\bar{\omega}^{ms} \bar{\omega}^{nj}\beta _{mn})\\
&=& \delta_s^k \bar{\omega}^{li}\bar{\omega}^{ms}\bar{\omega}^{nj}
\alpha _{kl}\beta _{mn}\\
&=& \bar{\omega}^{li}\bar{\omega}^{nj}
\bar{\omega}^{mk}\alpha _{kl}\beta _{mn}\\
&=& \bar{\omega}^{li}\bar{\omega}^{nj}
(-\bar{\omega}^{km})\alpha _{kl}\beta _{mn}\\
&=& - \bar{\omega}^{li}\bar{\omega}^{nj}(\alpha \diamond\beta)_{ln}
=\mu (\alpha \diamond \beta )^{ij}
\end{eqnarray*}

\medskip

\noindent 2.\\
\indent 2(a)
\begin{eqnarray*}
((\alpha\diamond\beta )\diamond\gamma )_{jl}&=&
\bar\omega^{r_2s_2}\bar\omega^{r_1s_1}
\alpha_{r_1r_2}\beta_{s_1j}\gamma_{s_2l}\\
&& \\
(\beta\diamond(\alpha\diamond\gamma ))_{jl}&=&
\bar\omega^{r_1s_1}\bar\omega^{r_2s_2}
\beta_{r_1j}\alpha_{r_2s_1}\gamma_{s_2l}\\
(r_1\leftrightarrow s_1) \ \ \ &=& \bar\omega^{s_1r_1}\bar\omega^{r_2s_2}
\beta_{s_1j}\alpha_{r_2r_1}\gamma_{s_2l}\\
&=& (-\bar\omega^{r_1s_1})\bar\omega^{r_2s_2}
\beta_{s_1j}(-\alpha_{r_1r_2})\gamma_{s_2l}=
((\alpha\diamond\beta)\diamond\gamma)_{jl}
\end{eqnarray*}

2(b)
$$(\tau\diamond\omega )_{ij}=\bar\omega^{rs}\tau_{ri}\omega_{sj}
=\delta_j^r\tau_{ri}=\tau_{ji}$$

\smallskip

Let us suppose, by induction, that, 
$\forall p \leq n$, $\alpha^{\diamond p}$ is
skewsymmetric and $\alpha^{\diamond n}
=\alpha^{\diamond l} \diamond \alpha^{\diamond m}$,
$\forall l,m \geq 1 \ s.t. \ l+m=n$.

\smallskip

\noindent We take now $l,m \geq 1 \ s.t. \ l+m=n+1$. We have, using 2(a), 2(b) 
and the induction hypothesis:
$$ \alpha^{\diamond (n+1)}=\alpha\diamond\alpha^{\diamond n}
=\alpha\diamond (\alpha^{\diamond (l-1)}\diamond\alpha^{\diamond m} )
=(\alpha^{\diamond (l-1)}\diamond\alpha )\diamond\alpha^{\diamond m}
=\alpha^{\diamond l}\diamond\alpha^{\diamond m} $$

and
\begin{eqnarray*}
(\alpha\diamond(\alpha^{\diamond n}\diamond\omega ))_{ij}
&=&(\alpha\diamond (-\alpha^{\diamond n}))_{ij} 
= - (\alpha^{\diamond (n+1)})_{ij}\\
(\alpha\diamond(\alpha^{\diamond n}\diamond\omega ))_{ij}
&=&((\alpha^{\diamond n}\diamond\alpha)\diamond\omega)_{ij}
=(\alpha^{\diamond (n+1)}\diamond\omega)=(\alpha^{\diamond (n+1)})_{ji}
\end{eqnarray*}

\hfill $\blacksquare $


\section{About Section 4}

\begin{lem}[Lemma 1]
$$\delta^{-1}\left(\frac{i}{2\hbar}[\alpha^{\diamond m}_{ij}y^i dx^j,
\alpha^{\diamond n}_{kl}y^kdx^l]\right)
=\frac{1}{2}\alpha^{\diamond (m+n)}_{il}y^i dx^l $$
\end{lem}

\noindent {\bf Proof:}
\begin{eqnarray*}
\frac{i}{2\hbar}[\alpha^{\diamond m}_{ij}y^i dx^j,
\alpha^{\diamond n}_{kl}y^kdx^l]&=&2\ \frac{i}{2\hbar}
(\alpha^{\diamond m}_{ij}y^i dx^j \circ_1 \alpha^{\diamond n}_{kl}y^kdx^l) \ \ \
\hbox{(Lemma A.1)}\\
&=&\frac{-i\hbar}{2}\ \frac{i}{\hbar}\bar\omega^{ik}\alpha^{\diamond m}_{ij}
\alpha^{\diamond n}_{kl}dx^j\wedge dx^l\\
&=& \frac{1}{2} \alpha^{\diamond (m+n)}_{il} dx^j\wedge dx^l
\end{eqnarray*}

\noindent and $\delta^{-1}\left(
\frac{1}{2} \alpha^{\diamond (m+n)}_{il} dx^j\wedge dx^l\right)
=\frac{1}{2} \alpha^{\diamond (m+n)}_{il} y^j dx^l$

\hfill $\blacksquare $

\begin{lem}[Lemma 2]
$$\delta^{-1}(\frac{i}{\hbar}[\alpha^{\diamond m}_{ij}y^i dx^j\
,\ \bar\omega^{kr}\alpha^{\diamond n}_{kl}\partial_rf\ y^l])
= \bar\omega^{kr}\alpha^{\diamond (m+n)}_{kj}\partial_rf\ y^j$$
\end{lem}
\noindent {\bf Proof:}
\begin{eqnarray*}
\delta^{-1}(\frac{i}{\hbar}[\alpha^{\diamond m}_{ij}y^i dx^j\
,\ \bar\omega^{kr}\alpha^{\diamond n}_{kl}\partial_rf\ y^l])
&=&2\ \frac{i}{\hbar}(\alpha^{\diamond m}_{ij}y^i dx^j
\circ_1 \bar\omega^{kr}\alpha^{\diamond n}_{kl}\partial_rf\ y^l)\\
&=&\frac{-i\hbar}{2}\ 2\frac{i}{\hbar} \bar\omega^{il} \bar\omega^{kr} 
\alpha^{\diamond m}_{ij} \alpha^{\diamond n}_{kl} \partial_rf dx^j\\
&=&\bar\omega^{kr} \bar\omega^{il} 
\alpha^{\diamond m}_{ij} (- \alpha^{\diamond n}_{lk}) \partial_rf dx^j\\
&=&\bar\omega^{kr}(- \alpha^{\diamond (m+n)}_{jk}) \partial_rf dx^j\\
&=&\bar\omega^{kr} \alpha^{\diamond (m+n)}_{kj} \partial_rf dx^j
\end{eqnarray*}
\noindent and $\delta^{-1}\left(
\bar\omega^{kr} \alpha^{\diamond (m+n)}_{kj} \partial_rf dx^j\right)
=\bar\omega^{kr} \alpha^{\diamond (m+n)}_{kj} \partial_rf y^j$

\hfill $\blacksquare $

\begin{lem}[Lemma 3]
$$\sigma\left( \bar\omega^{i_1l_1}
\alpha^{\diamond m}_{i_1j_1}\partial_{l_1}f\ y^{j_1}
\circ  \bar\omega^{i_2l_2}\alpha^{\diamond n}_{i_2j_2}
\partial_{l_2}g\ y^{j_2}\right)
= \frac{i \hbar}{2} (\bar{\alpha}^{\diamond(m+n)})^{l_1l_2} 
\partial_{l_1}f\partial_{l_2}g$$
\end{lem}

\noindent {\bf Proof:}
$\bar\omega^{i_1l_1}
\alpha^{\diamond m}_{i_1j_1}\partial_{l_1}f\ y^{j_1}
\circ  \bar\omega^{i_2l_2}\alpha^{\diamond n}_{i_2j_2}
\partial_{l_2}g\ y^{j_2}$ has a term given by $\circ_0$ which contains 
some $``y"$ and one given by $\circ_1$ which does not contain any. So
\begin{eqnarray*}
\sigma\big( \bar\omega^{i_1l_1}
\alpha^{\diamond m}_{i_1j_1}\partial_{l_1}f\ y^{j_1}
&\circ&\bar\omega^{i_2l_2}\alpha^{\diamond n}_{i_2j_2}
\partial_{l_2}g\ y^{j_2}\big)=\\
\hbox{(Lemma A.1)} \ \ \ &=&\bar\omega^{i_1l_1}
\alpha^{\diamond m}_{i_1j_1}\partial_{l_1}f\ y^{j_1}
\circ_1\bar\omega^{i_2l_2}\alpha^{\diamond n}_{i_2j_2}
\partial_{l_2}g\ y^{j_2}\\
&=&\frac{-i\hbar}{2}\ \bar\omega^{j_1j_2}\bar\omega^{i_1l_1}\bar\omega^{i_2l_2}
\alpha^{\diamond m}_{i_1j_1}\alpha^{\diamond n}_{i_2j_2}
\partial_{l_1}f\partial_{l_2}g\\
&=&\frac{-i\hbar}{2}\ \bar\omega^{i_1l_1}\bar\omega^{i_2l_2}
\alpha^{\diamond (m+n)}_{i_1i_2}\partial_{l_1}f\partial_{l_2}g\\
&=&\frac{i\hbar}{2}\ (- \bar\omega^{i_1l_1}\bar\omega^{i_2l_2}
\alpha^{\diamond (m+n)}_{i_1i_2})\partial_{l_1}f\partial_{l_2}g\\
&=&\frac{i\hbar}{2}\ 
(\bar\alpha^{\diamond (m+n)})^{l_1l_2}\partial_{l_1}f\partial_{l_2}g\\
\end{eqnarray*}

\hfill $\blacksquare $


\section{About Section 5}

We need the following notations:
$$R=\frac{1}{4} R_{ijkl}y^i y^j dx^k\wedge dx^l \ \
\hbox{so}\ \ \delta^{-1}R=\frac{1}{8} R_{ijkl}y^i y^j y^k dx^l$$
and 
\begin{equation} \label{24}
\delta^{-1}\frac{i}{\hbar} [\delta^{-1}R , a^{(1)} ]=
\delta^{-1}\frac{i}{\hbar} [\delta^{-1}R , y^m \partial_m f ]=
\frac{-1}{24}\bar\omega^{lm}R_{ijkl}y^i y^j y^k \partial_m f
\end{equation}
{\bf Proof of (\ref{24}):}
\begin{eqnarray*}
\delta^{-1}\frac{i}{\hbar} [\delta^{-1}R , a^{(1)} ]&=&
2\ \delta^{-1}R \circ_1 a^{(1)} \\
&=&\frac{-i\hbar}{2} 2 \frac{1}{8} 
(\bar\omega^{im}R_{ijkl}y^jy^k + \bar\omega^{jm}R_{ijkl}y^iy^k \\
&& \ \ +\bar\omega^{km}R_{ijkl}y^iy^j)\partial_mf dx^l\\
&=&\frac{-i\hbar}{8} \bar\omega^{im} (R_{ijkl} + R_{jikl} + R_{kjil})
y^jy^k\partial_mf dx^l\\
&=&\frac{-i\hbar}{8} \bar\omega^{im} (2R_{ijkl} + R_{kjil})
y^jy^k\partial_mf dx^l
\end{eqnarray*}
So 
$$ \delta^{-1}\frac{i}{\hbar} [\delta^{-1}R , a^{(1)} ]=
\frac{1}{3}\frac{1}{8} \bar\omega^{im} (2R_{ijkl} + R_{kjil})
y^jy^ky^l\partial_mf $$
By the symmetry properties of $R$ we have 
$R_{ijkl}=-R_{kjli}+R_{ljki}$ but $R_{kjli}y^jy^ky^l=-R_{ljki}y^jy^ky^l$
so $R_{ijkl}y^jy^ky^l=0$ and then
\begin{eqnarray*}
\delta^{-1}\frac{i}{\hbar} [\delta^{-1}R , a^{(1)} ]&=&
\frac{1}{24} \bar\omega^{im} R_{kjil}y^jy^ky^l\partial_mf \\
&=&\frac{1}{24} \bar\omega^{im} (-R_{kjli})y^jy^ky^l\partial_mf \\
&=&\frac{-1}{24} \bar\omega^{lm} R_{ijkl}y^iy^jy^k\partial_mf 
\end{eqnarray*}
by renaming the indices.

\hfill $\blacksquare $

\medskip

We define ${\cal{R}}_{l_1 l_2}$ by 
\begin{eqnarray*}
\sigma\left(\delta^{-1}R \circ \delta^{-1}R \right)
&=& \delta^{-1}R \circ_3 \delta^{-1}R \\
&=& \left(\frac{1}{8}\right)^2 {\cal{R}}_{l_1 l_2} dx^{l_1}\wedge dx^{l_2}
\end{eqnarray*}
i.e. $$ {\cal{R}}_{l_1l_2} = R_{i_1j_1k_1l_1} y^{i_1} y^{j_1} y^{k_1}\ 
\circ_3 \ R_{i_2j_2k_2l_2} y^{i_2} y^{j_2} y^{k_2}$$
In particular ${\cal{R}}_{l_1l_2}=-\ {\cal{R}}_{l_2l_1}$\ \ \ (Lemma A.1)

\medskip

\noindent We denote ${\cal{R}}^{m_1m_2} = -\bar{\omega}^{l_1m_1}
\bar{\omega}^{l_2m_2}{\cal{R}}_{l_1l_2}$

\bigskip

In Section 5 we defined $\beta_n=\frac{i}{\hbar^{n+3}}
\Big( (\delta^{-1}\partial)^n\delta^{-1}R 
\circ_{n+3} (\delta^{-1}\partial)^n\delta^{-1}R\Big)$.

\smallskip

\noindent Since $(\delta^{-1}\partial)^n\delta^{-1}R$
is a 1-form, for $n+3$ even ($n$ odd), we have, by Lemma A.1 :
\begin{eqnarray*}
(\delta^{-1}\partial)^n\delta^{-1}R 
& \circ_{n+3}& (\delta^{-1}\partial)^n\delta^{-1}R \\
&=& (-1)^{n+4} \ (\delta^{-1}\partial)^n\delta^{-1}R 
\ \circ_{n+3}\ (\delta^{-1}\partial)^n\delta^{-1}R \\
&=&\ - \ (\delta^{-1}\partial)^n\delta^{-1}R 
\ \circ_{n+3}\ (\delta^{-1}\partial)^n\delta^{-1}R
\end{eqnarray*}
so $\beta_n=0$ for $n$ odd.

\bigskip

We will now show:\\
(1)
$$\sigma\left(\delta^{-1}\frac{i}{\hbar}[\delta^{-1}R,a^{(1)}]\circ 
\delta^{-1}\frac{i}{\hbar}[\delta^{-1}R,b^{(1)}]\right)
= - \frac{1}{9.2^6} {\cal{R}}^{m_1m_2} \partial_{m_1}f \partial_{m_2}g$$

\medskip

\noindent (2) 
$$\begin{array}{ll} \sigma\Big(\delta^{-1}\frac{i}{\hbar}[\delta^{-1}R&,
\delta^{-1}\frac{i}{\hbar}[\delta^{-1}R,a^{(1)}]]\circ b^{(1)}\\
&+a^{(1)}\circ \delta^{-1}\frac{i}{\hbar}[\delta^{-1}R,
\delta^{-1}\frac{i}{\hbar}[\delta^{-1}R,b^{(1)}]]\Big)
\end{array}$$
$$= - \frac{1}{3.2^5} {\cal{R}}^{m_1m_2} 
\partial_{m_1}f \partial_{m_2}g $$

\medskip

\noindent (3)
$$\sigma\Big(\delta^{-1}\frac{i}{\hbar}[\hbar^{2}\delta^{-1}\beta_0,
a^{(1)}]\circ b^{(1)} + a^{(1)}\circ \delta^{-1}\frac{i}{\hbar}
[\hbar^{2}\delta^{-1}\beta_0,b^{(1)}]\Big)$$
$$= - \frac{1}{2^6} {\cal{R}}^{m_1m_2} \partial_{m_1}f \partial_{m_2}g$$

\medskip

{\bf Proofs:}\\
Proof of (1):
\begin{eqnarray*}
\sigma & \Big( & \delta^{-1}\frac{i}{\hbar}[\delta^{-1}R,a^{(1)}]\circ 
\delta^{-1}\frac{i}{\hbar}[\delta^{-1}R,b^{(1)}]\Big)\\
&=&(\frac{-1}{24})^2 
\bar\omega^{l_1m_1}R_{i_1j_1k_1l_1} 
y^{i_1} y^{j_1} y^{k_1} \partial_{m_1} f
\circ_3 
\bar\omega^{l_2m_2}R_{i_2j_2k_2l_2} 
y^{i_2} y^{j_2} y^{k_2} \partial_{m_2} g \\
&=& - \frac{1}{9.2^6} (-\bar{\omega}^{l_1m_1}
\bar{\omega}^{l_2m_2}{\cal{R}}_{l_1l_2}) \partial_{m_1} f \partial_{m_2} g \\
&=& - \frac{1}{9.2^6} {\cal{R}}^{m_1m_2} \partial_{m_1}f \partial_{m_2}g
\end{eqnarray*}

\medskip

\noindent Proof of (2): We have
\begin{eqnarray*}
\sigma\Big(\delta^{-1}\frac{i}{\hbar}[\delta^{-1}R&,&
\delta^{-1}\frac{i}{\hbar}[\delta^{-1}R,a^{(1)}]]\circ b^{(1)}\Big)\\
&=&\delta^{-1}\frac{i}{\hbar}[\delta^{-1}R,
\delta^{-1}\frac{i}{\hbar}[\delta^{-1}R,a^{(1)}]]_3\circ_1 b^{(1)}
\end{eqnarray*}
where
\begin{eqnarray*}
[\delta^{-1}R&,&
\delta^{-1}\frac{i}{\hbar}[\delta^{-1}R,a^{(1)}]]_3 \\
&=&\delta^{-1}R \circ_3
\delta^{-1}\frac{i}{\hbar}[\delta^{-1}R,a^{(1)}]
- \delta^{-1}\frac{i}{\hbar}[\delta^{-1}R,a^{(1)}] \circ_3 \delta^{-1}R \\
&=&2 \delta^{-1}R \circ_3
\delta^{-1}\frac{i}{\hbar}[\delta^{-1}R,a^{(1)}]\\
&=& 2 (\frac{-1}{24}) (\frac{1}{8}) 
R_{i_1j_1k_1l_1} 
y^{i_1} y^{j_1} y^{k_1} dx^{l_1}
\circ_3 
\bar\omega^{l_2m_2}R_{i_2j_2k_2l_2} 
y^{i_2} y^{j_2} y^{k_2} \partial_{m_2} f \\
&=&- \frac{1}{3.2^5} \bar\omega^{l_2m_2} 
{\cal{R}}_{l_1l_2} \partial_{m_2}f dx^{l_1}
\end{eqnarray*}
then
$$ \delta^{-1}\frac{i}{\hbar}[\delta^{-1}R,
\delta^{-1}\frac{i}{\hbar}[\delta^{-1}R,a^{(1)}]]_3
=- \frac{1}{3.2^5} \frac{i}{\hbar} \bar\omega^{l_2m_2} 
{\cal{R}}_{l_1l_2}\partial_{m_2}f  y^{l_1}$$
So,
\begin{eqnarray*}
\sigma\Big(\delta^{-1}\frac{i}{\hbar}[\delta^{-1}R&,&
\delta^{-1}\frac{i}{\hbar}[\delta^{-1}R,a^{(1)}]]\circ b^{(1)}\Big)\\
&=&- \frac{1}{3.2^5} \frac{i}{\hbar} \bar\omega^{l_2m_2} 
{\cal{R}}_{l_1l_2} \partial_{m_2}f y^{l_1} \circ_1 y^{m_1}\partial_{m_1}g\\
&=&\frac{-i\hbar}{2} (- \frac{1}{3.2^5}) \frac{i}{\hbar} 
\bar\omega^{l_1m_1}\bar\omega^{l_2m_2} 
{\cal{R}}_{l_1l_2} \partial_{m_2} f \partial_{m_1}g\\
(l_1\leftrightarrow l_2;m_1\leftrightarrow m_2)
&=&\frac{1}{2}(- \frac{1}{3.2^5})\bar\omega^{l_2m_2}\bar\omega^{l_1m_1}
{\cal{R}}_{l_2l_1}\partial_{m_1} f \partial_{m_2}g\\
&=&\frac{1}{2}(- \frac{1}{3.2^5})\bar\omega^{l_2m_2}\bar\omega^{l_1m_1}
(-{\cal{R}}_{l_1l_2})\partial_{m_1} f \partial_{m_2}g\\
&=&\frac{1}{2}(- \frac{1}{3.2^5}){\cal{R}}^{m_1m_2} 
\partial_{m_1}f \partial_{m_2}g
\end{eqnarray*}
In the same way
\begin{eqnarray*}
\sigma\Big(
a^{(1)}&\circ & \delta^{-1}\frac{i}{\hbar}[\delta^{-1}R,
\delta^{-1}\frac{i}{\hbar}[\delta^{-1}R,b^{(1)}]]\Big)\\
&=& \frac{1}{2}(- \frac{1}{3.2^5}) {\cal{R}}^{m_1m_2} 
\partial_{m_1}f \partial_{m_2}g  
\end{eqnarray*}
So we have (2).

\medskip

\noindent Proof of (3): We have $\hbar^2 \delta^{-1}\beta_0=\frac{1}{2^6}
\frac{i}{\hbar} {\cal R}_{l_1l_2} y^{l_1} dx^{l_2}$. Then
\begin{eqnarray*}
[\hbar^{2}\delta^{-1}\beta_0,a^{(1)}]
&=&2\ \hbar^{2}\delta^{-1}\beta_0 \circ_1 a^{(1)} \\
&=&(\frac{-i\hbar}{2}) 2 \frac{1}{2^6}
\frac{i}{\hbar} \bar\omega^{l_1m_1} {\cal{R}}_{l_1l_2} 
\partial_{m_1}f dx^{l_2}\\
&=&\frac{1}{2^6} \bar\omega^{l_1m_1} {\cal{R}}_{l_1l_2} 
\partial_{m_1}f dx^{l_2}\\
\hbox{so}\ \ \delta^{-1}\frac{i}{\hbar}[\hbar^{2}\delta^{-1}\beta_0,a^{(1)}]
&=& \frac{i}{\hbar}\frac{1}{2^6} \bar\omega^{l_1m_1} {\cal{R}}_{l_1l_2} 
\partial_{m_1}f y^{l_2}\ \ \hbox{and} \\ 
\sigma\Big(\delta^{-1}\frac{i}{\hbar}[\hbar^{2}\delta^{-1}\beta_0,
a^{(1)}]\circ b^{(1)}\Big)
&=&\frac{i}{\hbar}\frac{1}{2^6} \bar\omega^{l_1m_1} {\cal{R}}_{l_1l_2} 
\partial_{m_1}f y^{l_2}\circ_1 b^{(1)}\\
&=&(\frac{-i\hbar}{2})\frac{i}{\hbar}\frac{1}{2^6} 
\bar\omega^{l_2m_2}\bar\omega^{l_1m_1} {\cal{R}}_{l_1l_2} 
\partial_{m_1}f \partial_{m_2}g\\
&=&\frac{1}{2}(-\frac{1}{2^6} )
{\cal{R}}^{m_1m_2}\partial_{m_1}f \partial_{m_2}g
\end{eqnarray*}
In the same way
$$ \sigma\Big( a^{(1)}\circ \delta^{-1}\frac{i}{\hbar}
[\hbar^{2}\delta^{-1}\beta_0,b^{(1)}]\Big) =
\frac{1}{2}(-\frac{1}{2^6} )
{\cal{R}}^{m_1m_2}\partial_{m_1}f \partial_{m_2}g $$
So we have (3).

\hfill $\blacksquare $

\bigskip

\noindent {\bf About $\tilde{C}_{k+2}$:}

\medskip

The assertion is that there are no 1-differentiable terms in
$\tilde{C}_{k+2}$, for $k\geq 2$ 
(for $k=1$, Proposition \ref{p} shows that $\tilde{C}_3$ contains
$\frac{i}{2} \bar\alpha^{\diamond 2}$):
\begin{eqnarray*}
\tilde{C}_{k+2}=C_{k+2} + \sigma 
& \Big( &\underbrace{a^{(1)} \circ \rho^b_3 + \rho^a_3 \circ b^{(1)} }_{(i)}  \\
&+&\underbrace{a^{(2)} \circ \rho^b_2 + \rho^a_2 \circ b^{(2)} }_{(ii)}  \\
&+&\underbrace{a^{(3)} \circ \rho^b_1 
+ \rho^a_1 \circ b^{(3)} }_{(iii)}\ \Big) \\
\hbox{with}\ \ \ \ \ & & \\
\tilde{a}^{(2k+1)}
&=& a^{(2k+1)} + \rho^a_1\\
\hbox{so}\ \ \rho^a_1 &=& \frac{i}{\hbar}\delta^{-1}
[\hbar^k\delta^{-1}\alpha , a^{(1)}]\\
&& \\
\tilde{a}^{(2k+2)}
&=& a^{(2k+2)} + \rho^a_2\\
\hbox{so}\ \ \rho^a_2 &=& \delta^{-1}\partial \rho^a_1\\
&+&\frac{i}{\hbar}\delta^{-1}
\Big([\hbar^k\delta^{-1}\partial\delta^{-1}\alpha,a^{(1)}]
+ [\hbar^k\delta^{-1}\alpha , a^{(2)}] \Big)\\
&& \\
\tilde{a}^{(2k+3)}
&=& a^{(2k+3)} + \rho^a_3\\
\hbox{so}\ \ \rho^a_3 &=& \delta^{-1}\partial \rho^a_2\\
&+&\frac{i}{\hbar}\delta^{-1}
\Big( [\hbar^k\delta^{-1}\partial\delta^{-1}\partial
\delta^{-1}\alpha,a^{(1)}] +
[\hbar^k\delta^{-1}\partial\delta^{-1}\alpha,a^{(2)}]\\
&&\ \ +[\hbar^k\delta^{-1}\alpha , a^{(3)}] +
[\frac{i}{\hbar}\delta^{-1}[\delta^{-1}R,\hbar^k\delta^{-1}\alpha],a^{(1)}] \\
&&\ \ + 
[\delta^{-1}R,\frac{i}{\hbar}\delta^{-1}[\hbar^k\delta^{-1}\alpha,a^{(1)}]] 
\Big)
\end{eqnarray*}

\begin{itemize}
\item 
Since each term of $a^{(3)}$ contains three $y$'s and $\rho^a_1$ just one,
$\sigma$ yields (iii) to zero.

\item
$a^{(2)}= \frac{1}{2} y^iy^j\partial_i\partial_jf$ so $a^{(2)}$ 
and $b^{(2)}$ are 2-differentiable. Then (ii) is at least 2-differentiable
in one argument.

\item
Except a part of 
$\frac{i}{\hbar}\delta^{-1}
[\hbar^k\delta^{-1}\partial\delta^{-1}\alpha,a^{(2)}]$ which is 
2-differentiable, each term of $\rho^a_3$ contains three $y$'s.
Their product with $b^{(1)}$ contain two $y$'s. So they cancel applying 
$\sigma$. Thus (i) gives no 1-differentiable terms.
\end{itemize}

\end{document}